\documentclass[a4paper,12pt]{amsart}

\usepackage{amsmath}
\usepackage{amssymb}
\usepackage{mathrsfs}
\usepackage{ifthen}
\usepackage{graphicx}
\usepackage[T1]{fontenc} %skandit

\def\switchlinenumbers{\@ifstar
	{\let\makeLineNumberOdd\makeLineNumberRight
		\let\makeLineNumberEven\makeLineNumberLeft}%
	{\let\makeLineNumberOdd\makeLineNumberLeft
		\let\makeLineNumberEven\makeLineNumberRight}%
}

\def\setmakelinenumbers#1{\@ifstar
	{\let\makeLineNumberRunning#1%
		\let\makeLineNumberOdd#1%
		\let\makeLineNumberEven#1}%
	{\ifx\c@linenumber\c@runninglinenumber
		\let\makeLineNumberRunning#1%
		\else
		\let\makeLineNumberOdd#1%
		\let\makeLineNumberEven#1%
		\fi}%
}

\nonstopmode \numberwithin{equation}{section}
\newtheorem*{theorem*}{Theorem}

\newtheorem{thm}{Theorem}[section]
\newtheorem{cor}[equation]{Corollary}
\newtheorem{lem}[equation]{Lemma}

\theoremstyle{definition}

\newtheorem{prob}[equation]{Problem}
\newtheorem{rem}{Remark}[section]

\newtheorem{innercustomthm}{Theorem}
\newenvironment{customthm}[1]
{\renewcommand\theinnercustomthm{#1}\innercustomthm}
{\endinnercustomthm}

\newcounter{minutes}
\divide\time by 60
\newcounter{hours}
\multiply\time by 60
\addtocounter{minutes}{-\time}
\newcounter {own}
\def\theown {\thesection       .\arabic{own}}

\newenvironment{pf}[1][]{%
	\vskip 3mm
	\noindent
	\ifthenelse{\equal{#1}{}}%
	{{\slshape Proof. }}%
	{{\slshape #1.} }%
}%
{\qed\bigskip}

\newcounter{alphabet}

\def\be{\begin{equation}}
\def\ee{\end{equation}}

\newcommand{\bee}{\begin{enumerate}}
	\newcommand{\eee}{\end{enumerate}}

\newcommand{\blem}{\begin{lem}}
	\newcommand{\elem}{\end{lem}}
\newcommand{\bthm}{\begin{thm}}
	\newcommand{\ethm}{\end{thm}}
\newcommand{\bcor}{\begin{cor}}
	\newcommand{\ecor}{\end{cor}}
\newcommand{\beg}{\begin{examp}}
	\newcommand{\eeg}{\end{examp}}
\newcommand{\begs}{\begin{examples}}
	\newcommand{\eegs}{\end{examples}}
\newcommand{\bdefe}{\begin{defin}}
	\newcommand{\edefe}{\end{defin}}
\newcommand{\bprob}{\begin{prob}}
	\newcommand{\eprob}{\end{prob}}
\newcommand{\bei}{\begin{itemize}}
	\newcommand{\eei}{\end{itemize}}
\newcommand{\real}{{\operatorname{Re}\,}}

\newcommand{\norm}[1]{\left\lVert#1\right\rVert}
\newcommand{\abs}[1]{\left\lvert#1\right\rvert}

\begin{document}
	
	\title{Bohr operator on operator valued polyanalytic functions on simply connected domains}

	\author{Vasudevarao Allu}
	\address{Vasudevarao Allu,
		School of Basic Sciences,
		Indian Institute of Technology Bhubaneswar,
		Bhubaneswar-752050, Odisha, India.}
	\email{avrao@iitbbs.ac.in}
	
	\author{Himadri Halder}
	\address{Himadri Halder,
		School of Basic Sciences,
		Indian Institute of Technology Bhubaneswar,
		Bhubaneswar-752050, Odisha, India.}
	\email{himadrihalder119@gmail.com}

	\subjclass[{AMS} Subject Classification:]{Primary 47A56, 30B10, 47A63, ; Secondary 30C45}
	
	\def\thefootnote{}
	\footnotetext{ {\tiny File:~\jobname.tex,
			printed: \number\year-\number\month-\number\day,
			\thehours.\ifnum\theminutes<10{0}\fi\theminutes }
	} \makeatletter\def\thefootnote{\@arabic\c@footnote}\makeatother

	\begin{abstract}
	In this article, we study the Bohr operator for the operator valued subordination class $S(f)$ consisting of holomorphic functions subordinate to $f$ in the unit disk $\mathbb{D}:=\{z \in \mathbb{C}: |z|<1\}$, where $f:\mathbb{D} \rightarrow \mathcal{B}(\mathcal{H})$ is holomorphic and $\mathcal{B}(\mathcal{H})$ is the algebra of bounded linear operators on a complex Hilbert space $\mathcal{H}$. We establish several subordination results, which can be viewed as the analogues of a couple of interesting subordination results from scalar valued settings. We also obtain a von Neumann-type inequality for the class of self-analytic mappings of the unit disk $\mathbb{D}$ which fix the origin. Furthermore, we extensively study Bohr inequalities for operator valued polyanalytic functions in certain proper simply connected domains in $\mathbb{C}$. We obtain Bohr radius for the operator valued polyanalytic functions of the form $F(z)= \sum_{l=0}^{p-1} \overline{z}^l \, f_{l}(z) $, where $f_{0}$ is subordinate to an operator valued convex biholomorphic function, and operator valued starlike biholomorphic function in the unit disk $\mathbb{D}$.	
	\end{abstract}

	\maketitle
	\pagestyle{myheadings}
	\markboth{Vasudevarao Allu and  Himadri Halder}{Bohr operator on opertor valued polyanalytic functions on simply connected domains}
%\vspace{-6mm}
	
\section{Introduction}
Let $H^{\infty}(\mathbb{D})$ be the space of bounded analytic functions from the unit disk $\mathbb{D}:=\{z \in \mathbb{C}:|z|<1\}$ into the complex plane $\mathbb{C}$ and denote $\norm{f}_{\infty}:=\sup_{|z|<1} |f(z)|$. In $1914$, the following remarkable result for the universal constant $r=1/3$ for functions in $H^{\infty}(\mathbb{D})$ was proved by Harald Bohr \cite{Bohr-1914}.
\begin{customthm}{A} \label{him-P8-thm-A}
	Let $f \in H^{\infty}(\mathbb{D})$ with the power series $f(z)=\sum_{n=0}^{\infty} a_{n}z^{n}$. Then 
	\begin{equation} \label{him-p7-e-1.1}
	\sum_{n=0}^{\infty} |a_{n}|r^{n} \leq \norm{f}_{\infty}
	\end{equation}
	for $|z|=r \leq 1/3$ and the constant $1/3$, referred to as the classical Bohr radius, is the best possible.
\end{customthm}
%Initially, Bohr has proved this result for $r\leq 1/6$, which later has been improved independently to $1/3$ by M. Riesz, I. Schur, and F. Wiener. 
The interest in the Bohr inequality has been revived when Dixon \cite{Dixon & BLMS & 1995} used it to disprove the conjecture that if the non-unital von Neumann's inequality holds for a Banach algebra, then it is necessarily an operator algebra. In $2004$, Paulsen and Singh \cite{paulsen-2004-PAMS} extended Bohr's theorem to Banach algebras by finding a general version of Bohr inequality which is valid in the context of uniform algebras. For fixed $z\in \mathbb{D}$, we denote 
\begin{equation*}
\mathcal{G}_{z}:= \left\{f(z)=\sum_{n=0}^{\infty} a_{n}z^n: f\in H^{\infty}(\mathbb{D}) \right\}.
\end{equation*}
For $|z|=r$, the Bohr operator $M_{r}$ on $\mathcal{G}_{z}$ is defined by 
\begin{equation*}
M_{r}(f)= \sum_{n=0}^{\infty} |a_{n}| |z^n|= \sum_{n=0}^{\infty} |a_{n}| r^n.
\end{equation*}
The Bohr operator has the following properties, which has been established in \cite{paulsen-proof}.
\begin{thm} \label{him-P8-thm-1.1} \cite{abu-2021-JMMA}
For each fixed $z \in \mathbb{D}$ and $|z|=r$, the Bohr operator $M_{r}$ satisfies 
\begin{enumerate}
	\item $M_{r}(f) \geq 0$, and $M_{r}(f)=0$ if, and only if, $f\equiv 0$,
	\item  $M_{r}(f+g) \leq M_{r}(f)+M_{r}(g)$,
	\item $M_{r}(\alpha f)= |\alpha| M_{r}(f)$, $\alpha \in \mathbb{C}$,
	\item $M_{r}(f.g)\leq M_{r}(f).M_{r}(g)$,
	\item $M_{r}(1)=1$.
\end{enumerate}
\end{thm}
By the virtue of Theorem \ref{him-P8-thm-1.1}, it is worth to mention that the space $\mathcal{G}_{z}$ with the norm $M_{r}$ constitutes a Banach algebra. However, not all Banach spaces satisfy the Bohr phenomenon. In \cite{bene-2004}, B{\'e}n{\'e}teau {\it et al.} have shown that $H^{q}$, the usual Hardy spaces in $\mathbb{D}$ do not satisfy the Bohr phenomenon for any $0<q<\infty$.
% For fixed $z \in \mathbb{D}$, we denote the set 
%\begin{equation*}
%\mathcal{B}_{z}:= \{\psi \in \mathcal{G}_{z}: |\psi(z)|<1\}.
%\end{equation*}
%Then $\mathcal{B}_{z} \subset \mathcal{G}_{z}$ and $\mathcal{B}_{z}$ is also a Banach space. By virtue of Theorem \ref{him-P8-thm-A}, we observe that $\mathcal{B}_{z}$ is uniformly bounded when $|z| \leq 1/3$. 
A complex Banach algebra $\mathcal{A}$ satisfies the von Neumann inequality if for all polynomial $p(X)$ and for all $x \in \mathcal{A}$ with $\norm{x} \leq 1$, 
\begin{equation} \label{him-P8-e-1.2}
\norm{p(x)} \leq \norm{p}_{\infty}.
\end{equation}
In \cite{von-Neumann-1951}, von Neumann has shown that the algebra $\mathcal{L}(\mathcal{H})$ of all bounded operators on a Hilbert space $\mathcal{H}$ satisfies the inequality \eqref{him-P8-e-1.2}. It is well-known that every Banach algebra which is an operator algebra ({\it i.e.,} which is isometrically isomorphic to a closed subalgebra of $\mathcal{L}(\mathcal{H})$ for some Hilbert space $\mathcal{H}$) also satisfies the von Neumann inequality \eqref{him-P8-e-1.2}. It has been conjectured that 
\begin{enumerate}
	\item every Banach algebra is isomorphic to $\mathcal{L}(\mathcal{H})$,
	\item every Banach algebra satisfies the von Neumann inequality \eqref{him-P8-e-1.2}.
\end{enumerate}
In $1995$, Dixon \cite{Dixon & BLMS & 1995} disproved the conjecture (1). Bohr inequality has been extended to several complex variables and more abstract settings (see \cite{Ayt & Dja & BLMS & 2013,boas-1997,defant-2008,hamada-2009-Israel,paulsen-2002,paulsen-2004-PAMS,popescu-2019}).
%by considering the Banach algebra $l^{1}_{\beta}$ of sequences, 
%\begin{equation*}
%l^{1}_{\beta}= \left\{x=\{x_{n}\}^{\infty}_{n=1}: \frac{1}{\beta} \sum_{n=1}^{\infty} |x_{n}| < \infty\right\}.
%\end{equation*}
%Using Bohr theorem, Dixon has also shown that $l^{1}_{\beta}$ satisfies the non-unital von Neumann inequality if, and only if, $0<\beta\leq 1/3$.
\vspace{1.5mm}

Another interesting aspect of Bohr phenomenon thrives on considering the Bohr radius problem for subordinating families of analytic functions in $\mathbb{D}$. For two analytic functions $g$ and $f$ in $\mathbb{D}$, we say that $g$ is subordinate to $f$, written $g \prec f$, if there exists an analytic function $\phi:\mathbb{D} \rightarrow \mathbb{D}$ with $\phi(0)=0$ such that $g(z)=f(\phi(z))$ in $\mathbb{D}$. Let $S(f)$ be the class of analytic functions subordinate to $f$ in $\mathbb{D}$. We say that $g$ is quasi-subordinate to $f$ if there exists an analytic function $\psi$ with $|\psi(z)|\leq 1$ in $\mathbb{D}$ such that $g(z)=\psi(z)f(\phi(z))$ in $\mathbb{D}$. It is well-known that if $g$ is subordinate ( or quasi-subordinate) to $f$ in $\mathbb{D}$, then $M_{r}(g) \leq M_{r}(f)$ for $|z|=r \leq 1/3$. Bhowmik and Das \cite{bhowmik-2018} have studied the Bohr radius for the subordinating families and the Bohr radius for quasi-subordination families has been studied by Alkhaleefah {\it et al.} \cite{alkhaleefah-2019}. In $2021$, Bhowmik and Das \cite{bhowmik-2021} extended the Bohr phenomenon for the subordination to operator valued analytic functions in $\mathbb{D}$. %For discussing this, we need to introduce some basic notations and some definitions. 
\vspace{1mm}

Throughout this article, $\mathcal{B}(\mathcal{H})$ stands for the space of bounded linear operators on a complex Hilbert space $\mathcal{H}$. We want to concentrate operator valued holomoprhic functions $f:\mathbb{D} \rightarrow \mathcal{B}(\mathcal{H})$. The term subordination for operator valued functions can be defined as the scalar valued case. That is, for two holomorphic functions $g,f:\mathbb{D} \rightarrow \mathcal{B}(\mathcal{H})$, we say that $g$ is subordinate to $f$, written $g \prec f$, if there exists a holomorphic function $\phi:\mathbb{D} \rightarrow \mathbb{D}$ with $\phi(0)=0$ such that $g(z)=f(\phi(z))$ in $\mathbb{D}$. Let $S(f)$ be the class of analytic functions subordinate to $f$ in $\mathbb{D}$. For given two Banach spaces $X$ and $Y$ and a domain $\Omega \subset X$, a holomorphic function $f :\Omega \rightarrow Y$ is said to be biholomorphic on $\Omega$ if $f^{-1}$ exists and is holomorphic in $f(\Omega) \subseteq Y$. We say that a biholomorphic function $f$ is starlike in its domain $\Omega$ with respect to $\xi_{0} \in \Omega$ if $f(\Omega)$ is a starlike domain with respect to $f(\xi_{0})$ {\it i.e.,} $(1-t) f(\xi_{0}) + tf(z) \in f(\Omega)$ for all $z \in \Omega$ and $t \in [0,1]$. A biholomorphic function $f$ is called starlike if $f$ is starlike with respect to $0\in \Omega$ and $f(0)=0$. A biholomorphic function $f$ is said to be convex if $f$ is starlike with respect to every point in $\Omega$. For convex or starlike biholomorphic function $f$ in $\mathbb{D}$, Bohr phenomenon for any $g \in S(f)$ has been studied in \cite{bhowmik-2021}. 
\par
For the rest of our discussion, we introduce some notations. Throughout this paper, $\norm{A}$ stands for the operator norm of $A$ for any $A\in \mathcal{B}(\mathcal{H})$ and $\abs{A}=(A^{*}A)^{1/2}$ denotes the absolute value of $A$, where $A^{*}$ is the adjoint of $A$ and $T^{1/2}$ denotes the unique positive square root of a positive operator $T$. We denote $I$ be the identity operator on $\mathcal{H}$. 
\vspace{1.5mm}

%Bohr's classical theorem and its generalizations on simply connected domains are now active areas of research in numerous analytic function spaces. 
In $2010$, Fournier and Ruscheweyh \cite{Four-Rusc-2010} extensively studied the Bohr radius problem for arbitrary simply connected domains containing $\mathbb{D}$. Let $\mathcal{H}(\Omega)$ be the class of analytic functions $f : \Omega \rightarrow \mathbb{C}$ and $\mathcal{B}(\Omega)$ denote the class of functions $f \in \mathcal{H}(\Omega)$ such that $f(\Omega) \subseteq \overline{\mathbb{D}}$. For the class $\mathcal{B}(\Omega)$, the Bohr radius $\mathcal{B}_{\Omega}$ is defined by (see \cite{Ahamed-Allu-Halder-P3-2020,Four-Rusc-2010})
$$
B_{\Omega}:=\sup\bigg\{r\in (0,1) : M_{r}(f)\leq 1\; \text{for all}\; f(z)=\sum_{n=0}^{\infty}a_nz^n\in\mathcal{B}(\Omega),\; z\in\mathbb{D}\bigg\},
$$
where $M_{r}(f):=\sum_{n=0}^{\infty}|a_n|r^n$ is the Bohr operator for $f \in \mathcal{B}(\Omega)$ in $\mathbb{D}$. For $\Omega=\mathbb{D}$, $\mathcal{B}(\Omega)$ reduces to $B_{\mathbb{D}}=1/3$, which is the classical Bohr radius for the class $\mathcal{B}(\mathbb{D})$.
\vspace{2mm}

For $0\leq \gamma <1$, Fournier and Ruscheweyh \cite{Four-Rusc-2010} have estimated the Bohr radius for the class $\mathcal{B}(\Omega_{\gamma})$ and proved that $\mathcal{B}_{\Omega_{\gamma}}=(1+\gamma)/(3+\gamma)$, where 
$$
\Omega_{\gamma}:=\bigg\{z\in\mathbb{C} : \bigg|z+\frac{\gamma}{1-\gamma}\bigg|<\frac{1}{1-\gamma}\bigg\}.
$$
Let $H^{\infty}(\Omega,X)$ be the space of bounded analytic functions from $\Omega$ into a complex Banach space $X$ and $\norm{f}_{H^{\infty}(\Omega,X)}= \sup_{z \in \Omega} \norm{f(z)}$. 
The Bohr phenomenon for operator valued functions on simply connected domains  has been studied in \cite{Himadri-Vasu-P7}. Let $\mathcal{B}(\mathcal{H})$ be the algebra of all bounded linear operators on a complex Hilbert space $\mathcal{H}$. For the class $H^{\infty}(\Omega,\mathcal{B}(\mathcal{H}))$, we denote (see \cite{Himadri-Vasu-P7})
\begin{equation} \label{him-P7-e-1.24}
\lambda_{\mathcal{H}}:=\lambda_{\mathcal{H}}(\Omega):= \sup \limits _{\substack{f \in  H^{\infty}(\Omega,\mathcal{B}(\mathcal{H}))\\{\norm{f(z)}\leq1}}} \left\{\frac{\norm{A_{n}}}{\norm{\, I- |A_{0}|^2}\,} : A_{0} \not \equiv f(z)= \sum \limits_{n=0}^{\infty} A_{n} z^n, \,\, z \in \mathbb{D}\right\}.
\end{equation}
%The following results have been obtained in \cite{Himadri-Vasu-P7}.
Recently, Allu and Halder \cite{Himadri-Vasu-P7} have established the Bohr theorem for the functions in $H^{\infty}(\Omega,\mathcal{B}(\mathcal{H}))$.
\begin{thm} \label{him-P7-cor-1.40} \cite{Himadri-Vasu-P7}
	Let $f \in H^{\infty}(\Omega,\mathcal{B}(\mathcal{H}))$ with $\norm{f(z)}_{H^{\infty}(\Omega,\mathcal{B}(\mathcal{H}))} \leq 1$ such that $f(z)= \sum_{n=0}^{\infty} A_{n}z^n$ in $\mathbb{D}$, where $A_{0}=\alpha_{0}I$ for $|\alpha_{0}|<1$ and $A_{n} \in \mathcal{B}(\mathcal{H})$ for all $n \in \mathbb{N} \cup \{0\}$. Then 
	\begin{equation}
	\sum_{n=0}^{\infty} \norm{A_{n}} r^n \leq 1 \,\,\,\,\,\, \mbox{for} \,\,\,\, r \leq \frac{1}{1+2\lambda_{\mathcal{H}}}.
	\end{equation}
\end{thm}
For $\Omega=\Omega_{\gamma}$ and $p=1$ in \cite[Corollary 1.52]{Himadri-Vasu-P7}, we have the following result.
\begin{thm} \label{him-P7-cor-1.52} \cite{Himadri-Vasu-P7}
	Let $f \in H^{\infty}(\Omega_{\gamma},\mathcal{B}(\mathcal{H}))$ with $\norm{f(z)}_{H^{\infty}(\Omega_{\gamma},\mathcal{B}(\mathcal{H}))} \leq 1$ such that $f(z)= \sum_{n=0}^{\infty} A_{n}z^n$ in $\mathbb{D}$, where $A_{0}=\alpha_{0}I$ for $|\alpha_{0}|<1$ and $A_{n} \in \mathcal{B}(\mathcal{H})$ for all $n \in \mathbb{N} \cup \{0\}$. Then 
	\begin{equation}
	\sum_{n=0}^{\infty} \norm{A_{n}} r^n \leq 1 \,\,\,\,\,\, \mbox{for} \,\,\,\, r \leq \frac{1+\gamma}{3+\gamma}.
	\end{equation}
\end{thm}
When $\Omega_{\gamma}=\mathbb{D}$ {\it i.e.,} $\gamma=0$, under the same assumptions as in Theorem \ref{him-P7-cor-1.52}, we have 
\begin{equation} \label{him-P8-e-1.6}
\sum_{n=0}^{\infty} \norm{A_{n}} r^n \leq 1 \,\,\,\,\,\, \mbox{for} \,\,\,\, r \leq \frac{1}{3}.
\end{equation}
\vspace{-6mm}

\section{Bohr operator on operator valued subordination classes}
In this section, we study subordination results for Bohr operator on operator valued analytic functions in $\mathbb{D}$. Recall that $\mathcal{B}(\mathcal{H})$ be the algebra of all bounded linear operators on a complex Hilbert space $\mathcal{H}$. For analytic functions $f:\mathbb{D} \rightarrow \mathcal{B}(\mathcal{H})$ with $f(z)= \sum_{n=0}^{\infty} A_{n} z^n$ in $\mathbb{D}$ and $A_{n} \in \mathcal{B}(\mathcal{H})$ for $n \in \mathbb{N} \cup \{0\}$, we define the Bohr operator $M_{r}(f)$ as the scalar valued case. That is, $M_{r}(f)= \sum_{n=0}^{\infty} \norm{A_{n}}|z|^n$. It can be easily seen that the operator $M_{r}$ satisfies the same property as in Theorem \ref{him-P8-thm-1.1}. Infact, for $f,g :\mathbb{D} \rightarrow \mathcal{B}(\mathcal{H})$ with $f(z)= \sum_{n=0}^{\infty} A_{n} z^n$ and $g(z)= \sum_{n=0}^{\infty} B_{n} z^n$ in $\mathbb{D}$ with $A_{n}, B_{n} \in \mathcal{B}(\mathcal{H})$ for $n \in \mathbb{N} \cup \{0\}$, we have 
\begin{equation} \label{him-P8-e-2.1}
M_{r}(f+g) = \sum_{n=0}^{\infty} \norm{A_{n}+B_{n}} r^n \leq \sum_{n=0}^{\infty} \norm{A_{n}} r^n + \sum_{n=0}^{\infty} \norm{B_{n}} r^n= M_{r}(f) + M_{r}(g).
\end{equation}
Using \eqref{him-P8-e-2.1}, it is easy to see that if $F(z)= \sum_{k\in \mathbb{Z}} f_{k}(z)$ is analytic function in $\mathbb{D}$, then $M_{r}(F) \leq \sum_{k\in \mathbb{Z}} M_{r}(f_{k})$, where $f_{k}: \mathbb{D} \rightarrow \mathcal{B}(\mathcal{H})$ is analytic function in $\mathbb{D}$ for each $k \in \mathbb{Z}$. On the other hand, we note that $M_{r}(\beta f)= |\beta|M_{r}(f)$ for any $\beta \in \mathbb{C}$ and $M_{r}(z^pf)=r^p M_{r}(f)$. We observe that $(fg)(z)=\sum_{n=0}^{\infty} A_{n} (z^n g(z))$ and hence
\begin{equation} \label{him-P8-e-2.1-a}
M_{r}(fg) \leq \sum_{n=0}^{\infty} \norm{A_{n}} r^n M_{r}(g)= M_{r}(f)M_{r}(g).
\end{equation} 
Clearly, $M_{r}(I)=1$. The following result has been established in \cite{abu-2021-JMMA}.
\vspace{-1mm}

\begin{lem} \label{him-P8-lem-2.1} \cite{abu-2021-JMMA}
Let $\phi:\mathbb{D} \rightarrow \mathbb{D}$ be analytic function with $\phi(0)=0$. Then $M_{r}(\phi) \leq |z|$ for $|z|=r \leq 1/3$.
\end{lem}
%\begin{pf}
%Since $\phi:\mathbb{D} \rightarrow \mathbb{D}$ be analytic function with $\phi(0)=0$, by virtue of the Schwarz lemma, $\psi(z)=\phi(z)/z$ is again an analytic self map of $\mathbb{D}$. Then, the classical Bohr theorem gives that $M_{r}(\psi(z))=M_{r}(\phi(z)/z) \leq 1$ for $r\leq 1/3$ {\it i.e.,} $M_{r}(\phi) \leq |z|$ for $|z|=r \leq 1/3$.  
%\end{pf}
 The following result is the operator-valued subordination result for Bohr operator, which has been first proved in \cite{bhowmik-2021}. By using Lemma \ref{him-P8-lem-2.1}, we give an alternative proof.
\begin{thm} \label{him-P8-thm-2.1}
Let $f,g:\mathbb{D} \rightarrow \mathcal{B}(\mathcal{H})$ be holomorphic functions such that $f \prec g$. Then 
\begin{equation} \label{him-P8-e-2.2}
M_{r}(f) \leq M_{r}(g)\,\,\,\,\,\,\,\,\, \mbox{for} \,\,\, |z|=r\leq \frac{1}{3}.
\end{equation}
\end{thm}
\begin{pf}
Let $f(z)= \sum_{n=0}^{\infty} A_{n} z^n$ and $g(z)= \sum_{n=0}^{\infty} B_{n} z^n$ in $\mathbb{D}$ with $A_{n}, B_{n} \in \mathcal{B}(\mathcal{H})$ for $n \in \mathbb{N} \cup \{0\}$.
Since $f \prec g$ in $\mathbb{D}$, then there exists an analytic function $\phi:\mathbb{D} \rightarrow \mathbb{D}$ such that $\phi(0)=0$ and $f(z)=g(\phi(z))$ in $\mathbb{D}$. In view of \eqref{him-P8-e-2.1}, \eqref{him-P8-e-2.1-a}, and Lemma \ref{him-P8-lem-2.1}, for $0\leq |z|=r\leq 1/3$, we obtain
\begin{align*}
M_{r}(f)=M_{r}(g(\phi))&=M_{r}\left(\sum_{n=0}^{\infty}B_{n}(\phi(z))^n\right) \\&  
\leq \sum_{n=0}^{\infty} \norm{B_{n}}(M_{r}(\phi(z)))^n \leq \sum_{n=0}^{\infty} \norm{B_{n}}r^n=M_{r}(g).
\end{align*}
This completes the proof.
\end{pf}

In particular for the scalar valued functions  $f,g:\mathbb{D} \rightarrow \mathbb{C}$, Theorem \ref{him-P8-thm-2.1} reduces to the result of Abu Muhanna {\it et al.} \cite{abu-2021-JMMA}, and Bhowmik and Das \cite{bhowmik-2018}. In view of Theorem \ref{him-P8-thm-2.1}, we obtain the following interesting result.

\begin{thm} \label{him-P8-thm-2.2}
Let $f,g,h:\mathbb{D} \rightarrow \mathcal{B}(\mathcal{H})$ be holomorphic functions such that $f(z)=h(z)g(\phi(z))$ for some analytic function $\phi:\mathbb{D} \rightarrow \mathbb{D}$ with $\phi(0)=0$. If $\norm{h(z)} \leq M$ for $|z|<\beta \leq 1$ and $h(0)=\alpha I$ with $|\alpha| \leq M$, then $M_{r}(f) \leq M \, M_{r}(g)$ for $0 \leq r \leq \beta/3$.
\end{thm}  
\begin{pf}
From \eqref{him-P8-e-2.1-a}, we have 
\begin{equation} \label{him-P8-e-2.3}
M_{r}(f) \leq M_{r}(h) M_{r}(g(\phi)).
\end{equation}
The assumption $\norm{h(z)} \leq M$ in the disk $\mathbb{D}_{\beta}:=\{z \in \mathbb{C}:|z|<\beta\}$ shows that $h_{1}:\mathbb{D} \rightarrow \mathcal{B}(\mathcal{H})$ defined by  $h_{1}(z)=h(z)/M$ is holomorphic and $\norm{h_{1}(z)} \leq 1$ in $\mathbb{D}_{\beta}$ such that $h_{1}(0)=(\alpha /M)I$. Since $|\alpha| \leq M$, from \eqref{him-P8-e-1.6}, we obtain 
\begin{equation} \label{him-P8-e-2.4}
M_{r}(h) \leq M \,\,\,\,\,\mbox{for} \,\,\,\, 0<r\leq \frac{\beta}{3}.
\end{equation}
Furthermore, in view of Theorem \ref{him-P8-thm-2.1}, we have
\begin{equation} \label{him-P8-e-2.5}
M_{r}(g(\phi)) \leq M_{r}(g) \,\,\,\,\,\mbox{for} \,\,\,\, 0<r\leq\frac{1}{3}.
\end{equation}
By using \eqref{him-P8-e-2.4} and \eqref{him-P8-e-2.5} in \eqref{him-P8-e-2.3}, we obtain
\begin{equation}
M_{r}(f) \leq M\,\, M_{r}(g) \,\,\,\,\,\mbox{for} \,\,\,\, 0<r\leq \frac{\beta}{3}.
\end{equation}
This completes the proof.
\end{pf}
\vspace{-5mm}

\begin{rem}
\begin{enumerate}
	\item For a particular case $h(z)\equiv I$, Theorem \ref{him-P8-thm-2.2} reduces to Theorem \ref{him-P8-thm-2.1}. By taking $f,g,h :\mathbb{D} \rightarrow \mathbb{C}$ are analytic functions in Theorem \ref{him-P8-thm-2.2}, we obtain the scalar valued quasi-subordination result which has been established in \cite{alkhaleefah-2019}. 
	
	\item When $\norm{h(z)} \leq 1$ in $\mathbb{D}$, we deduce that 
	$M_{r}(f) \leq M_{r}(g)$ for $|z|=r \leq 1/3$.
\end{enumerate}
\end{rem}
We now prove the following interesting result, which is an analogue of von Neumann inequality \eqref{him-P8-e-1.2}.
\begin{thm} \label{him-P8-thm-2.3}
Let $f:\mathbb{D} \rightarrow \mathcal{B}(\mathcal{H})$ be analytic in $\mathbb{D}$ and continuous in $\overline{\mathbb{D}}$ such that $f(0)=\alpha I$ for some $\alpha \in \mathbb{C}$ with $|\alpha|<1$. Then 
\begin{equation} \label{him-P8-e-2.8}
M_{r}(f(\phi)) \leq \norm{f}_{\infty} \,\,\,\, \mbox{for} \,\,\,\,\, 0\leq r \leq 1/3,
\end{equation} 
where $\phi: \mathbb{D} \rightarrow \mathbb{D}$ is analytic function with $\phi(0)=0$.
\end{thm}
\begin{pf}
Let $f(z)= \sum_{n=0}^{\infty} A_{n} z^n$ in $\mathbb{D}$, where $A_{0}=\alpha I$ and $A_{n} \in \mathcal{B}(\mathcal{H})$ for $n \in \mathbb{N} \cup \{0\}$. Then, for $r \leq 1/3$, Theorem \ref{him-P8-thm-2.1} gives 
\begin{equation} \label{him-P8-e-2.9}
M_{r}(f(\phi)) \leq  \sum_{n=0}^{\infty} \norm{A_{n}} r^n= M_{r}(f).
\end{equation}
In view of \eqref{him-P8-e-1.6}, for $0\leq r \leq 1/3$, we obtain $M_{r}(f) \leq \norm{f}_{\infty}$ which together with \eqref{him-P8-e-2.9} gives \eqref{him-P8-e-2.8}.
\end{pf}
\vspace{-8mm}

\section{Bohr theorem for operator-valued polyanalytic functions }
Polyanalytic functions $f$ of order $p$ defined in a simply connected domain $\Omega \subseteq \mathbb{C}$ are complex-valued polynomials in the variable $\overline{z}$ with analytic functions are their coefficients. That is, $f$ has the following form 
$
f(z)= \sum_{l=0}^{p-1} \overline{z}^p \, f_{l}(z),
$
where $f_{l}$'s are analytic functions in $\Omega$. Equivalently, polyanalytic functions can also be defined as the $\mathcal{C}^p(\Omega)$-solutions of the generalized Cauchy-Riemann equations $\partial ^{p} f/ \partial \overline{z}^p=0$ in $\Omega$ (the Cauchy-Riemann equations of order $p$). Throughout this paper, we assume that $p \geq 2$. 
\vspace{2mm}

In $1908$, Kolossov \cite{Kolossov-1908} first introduced polyanalytic functions in connection with his research in the mathematical theory of elasticity. Polyanalytic function theory has been extensively studied by Balk \cite{balk}. In $2011$, Agranovsky \cite{agranovsky-2011} characterized the polyanalytic functions by meromorphic extensions into chains of circles. It is worth mentioning that the properties of polyanalytic functions can be different from those of analytic functions (see \cite{balk}). By considering the polyanalytic function $f(z)=1-z\overline{z}$, one can easily see that $f$ vanishes on the boundary of the unit disk $\mathbb{D}$ without vanishing identically in $\mathbb{D}$. Studying polyanalytic functions also reveals some new properties of analytic functions. The study of polyanalytic functions is closely related to numerous research topics of complex analysis {\it e.g.,} function theory of several complex variables, the theory of distribution of values of meromorphc functions, the theory of meromorphic curves, the theory of boundary properties of analytic functions. In $2019$, Hachadi and Youssfi \cite{Hachadi-CAOT-polyanalytic} have studied several properties of polyanalytic reproducing kernels. In $2021$, Abdulhadi and Hajj \cite{Abdulhadi-2021} extensively studied univalency criteria, Landau's theorem, arc-length problem, and the Bohr phenomenon problem for polyanalytic functions in $\mathbb{D}$.   
\vspace{2mm}

%In particular, a complex-valued $\mathcal{C}^2 (\Omega)$ in $\Omega$ is said to be bi-analytic if 
%\begin{equation} \label{him-P7-e-1.66}
%\frac{\partial}{\partial \overline{z}} \left(\frac{\partial f}{\partial \overline{z}}\right)=0.
%\end{equation}
%Then $f$ is a polyanalytic function of order $2$. From \eqref{him-P7-e-1.66}, we note that $f_{\overline{z}}$ is analytic in $\Omega$. It can be shown that $f$ has the form $f(z)=\overline{z}f_{1}(z) + f_{0}(z)$, where $f_{0}$ and $f_{1}$ are analytic functions in $\Omega$ and $f_{1}=f_{\overline{z}}$. 
Since complex-valued polyanalytic functions are polynomials in $\overline{z}$ in simply connected domain $\Omega$, this leads to study the the operator-valued polyanalytic functions. A operator-valued polyanalytic function $F$ of order $p$ in $\Omega$ is a polynomial in $\overline{z}$ with operator valued analytic function as its coefficients. That is, $F$ has the following form 
\begin{equation*}
F(z)= \sum_{l=0}^{p-1} \overline{z}^p \, f_{l}(z),
\end{equation*}
where $f_{l}: \Omega \rightarrow \mathcal{B}(\mathcal{H})$ are analytic functions for $l=0,1, \ldots , p-1$ and $f_{p-1} \not \equiv 0$. Now we consider the simply connected domain $\Omega$ containing $\mathbb{D}$.
\vspace{2mm}

Although, Bohr radius and Bohr phenomenon have been extensively studied, no attempt has been made so far to obtain operator valued analogues of Bohr phenomenon for polyanalytic functions. Therefore, our main aim of this section is to obtain the Bohr inequality under appropriate considerations and necessary conditions. In the following result we establish operator valued analogues of Bohr inequality in simply connected domain $\Omega$ containing $\mathbb{D}$.
\vspace{-2mm}

\begin{thm} \label{him-P7-thm-1.7}
	Let $F$ be a polyanalytic function of order $p$ in $\Omega$ with $F(z)= \sum_{l=0}^{p-1} \overline{z}^l \, f_{l}(z) $, where each $f_{l}: \Omega \rightarrow \mathcal{B}(\mathcal{H})$ is analytic function such that $f_{l}(z)= \sum_{n=0}^{\infty} A_{n,l} z^n$ in $\mathbb{D}$ and $A_{n,l} \in \mathcal{B}(\mathcal{H})$ for $n \in \mathbb{N} \cup \{0\}$. Also assume that 
	\begin{enumerate}
		\item [\rm{(a)}]$\norm{f_{0}(z)} \leq 1$ in $\Omega$ such that $f_{0}(0)={\bf 0}$ and $f'_{l}(0)=\alpha_{l} f'_{0}(0)$ with $|\alpha_{l}| <k$ for $k \in [0,1]$ and each $l=1,\ldots, p-1$.
		\item [\rm{(b)}] $\omega_{l} : \Omega \rightarrow \mathcal{B}(\mathcal{H})$ is analytic with $\norm{\omega_{l}(z)} \leq k$ in $\Omega$ for $k \in [0,1]$, where $\omega_{l}(z)= f'_{l}(z) (f'_{0}(z))^{-1}$ in $\Omega$ such that $\omega_{l}(z)=\sum_{n=0}^{\infty} \omega_{n,l} z^n $ in $\mathbb{D}$, provided $(f'_{0}(z))^{-1}$ exists for all  $z \in \Omega$.
	\end{enumerate}
	Then $M_{r}(F) \leq 1$ for $|z|=r \leq R_{f}= \min \{r_{f}(p), 1/(1+2\lambda _{\mathcal{H}})\}$, where $r_{f}(p)$ is the smallest root in $(0,1)$ of 
	\begin{equation} \label{him-P7-e-1.67}
	(1-r)^2 - \lambda_{\mathcal{H}}\, r - \lambda_{\mathcal{H}}\, r^{p+1}=0.
	\end{equation}
	Here $\lambda _{\mathcal{H}}$ is given by \eqref{him-P7-e-1.24}.
\end{thm}
\begin{pf}
	Let $F(z)= \sum_{l=0}^{p-1} \overline{z}^l \, f_{l}(z) $ with $f_{l}(z)= \sum_{n=0}^{\infty} A_{n,l} \,z^n$ in $\mathbb{D}$. Then 
	\begin{equation} \label{him-P7-e-1.68}
	M_{r}(F)= M_{r} \left(\sum_{l=0}^{p-1} \overline{z}^l \, f_{l}(z)\right) \leq \sum_{l=0}^{p-1} r^l M_{r}(f_{l}) \,\,\,\,\, \mbox{for} \,\,\,\, |z|=r<1.
	\end{equation}
	Since $\omega_{l}(z)=f'_{l}(z) (f'_{0}(z))^{-1}$ in $\Omega$ with  $\norm{\omega_{l}(z)} \leq k$ in $\Omega$ for each $l$ such that $f'_{l}(0)=\alpha_{l} f'_{0}(0)$, it follows that $f'_{l}(z)= \omega_{l}(z) f'_{0}(z)$ in $\Omega$ with $\omega_{l}(0)=\alpha_{l}I$, where $|\alpha_{l}|<k$ for each $l=1, \ldots , p-1$. Let $\lambda _{\mathcal{H}}$ be given by \eqref{him-P7-e-1.24}. In view of Theorem \ref{him-P7-cor-1.40}, for $|z|=r \leq 1/(1+2\lambda_{\mathcal{H}})$, we have $M_{r}(\omega_{l}) \leq k$, which together with \eqref{him-P8-e-2.1-a} gives 
	\begin{equation} \label{him-P7-e-1.69}
	M_{r}(f_{l})= \int \limits_{0}^{r} M_{r}(f'_{l})dt= \int\limits_{0}^{r} M_{r}(\omega_{l}f'_{0})dt \leq k\,\int\limits_{0}^{r} M_{r}(f'_{0})dt =k\, M_{r}(f_{0}).
	\end{equation}
	Using \eqref{him-P7-e-1.68} and \eqref{him-P7-e-1.69}, for $|z|=r \leq 1/(1+2\lambda_{\mathcal{H}})$, we obtain
	\begin{equation} \label{him-P7-e-1.70}
	M_{r}(F) \leq k\sum_{l=0}^{p-1} r^l M_{r}(f_{0}) = k M_{r}(f_{0}) \left(\frac{1-r^p}{1-r}\right).
	\end{equation}	
We now wish to find the upper bound for $M_{r}(f_{0})$. We observe that $f_{0}: \Omega \rightarrow \mathcal{B}(\mathcal{H})$ is analytic function with $\norm{f_{0}(z)}\leq 1$ in $\Omega$ such that $f_{0}(z)= \sum_{n=0}^{\infty} A_{n,0} z^n$ in $\mathbb{D}$, where $f_{0}(0)=A_{0,0}={\bf 0}$. Then in view of \eqref{him-P7-e-1.24}, we have $\norm{A_{n,0}} \leq \lambda_{\mathcal{H}}$ for $n \geq 1$ and hence 
	\begin{equation} \label{him-P7-e-1.71}
	M_{r}(f_{0}) = \sum_{n=0}^{\infty} \norm{A_{n,0}}\, r^n \leq \lambda_{\mathcal{H}} \left(\frac{r}{1-r}\right).
	\end{equation}
	In view of \eqref{him-P7-e-1.70} and \eqref{him-P7-e-1.71}, for $r \leq 1/(1+2\lambda_{\mathcal{H}})$, we obtain 
	\begin{equation} \label{him-P7-e-1.72}
	M_{r}(F)  \leq k \, \lambda_{\mathcal{H}}\left( \frac{r}{1-r}\right) \left(\frac{1-r^p}{1-r}\right).
	\end{equation}
	Therefore, $M_{r}(F) \leq 1$ for $r \leq \min \{1/(1+2\lambda_{\mathcal{H}}), r_{f}(p)\}$, where $r_{f}(p)$ is the smallest root in $(0,1)$ of
	$$
k\,	\lambda_{\mathcal{H}}\left( \frac{r}{1-r}\right) \left(\frac{1-r^p}{1-r}\right) = 1,
	$$
	which is equivalent to $(1-r)^2 -k\, \lambda_{\mathcal{H}}\,r + k\,\lambda_{\mathcal{H}}\, r^{p+1}=0$.
%	The existence of the root can be shown by considering the function $H:[0,1] \rightarrow \mathbb{R}$ defined by $H(r)=(1-r)^2 - k\lambda_{\mathcal{H}}\,r + k \lambda_{\mathcal{H}}\, r^{p+1}$. We note that $H$ is continuous in $[0,1]$ such that $H(0)=1>0$ and $H(1)=-2k\lambda_{\mathcal{H}}<0$. Then by the intermediate value theorem, $H$ has a root in $(0,1)$ and choose $r_{f}(p)$ to be the smallest root. 
This completes the proof.
\end{pf}

As a consequence of Theorem \ref{him-P7-thm-1.7}, we obtain Bohr type inequality for bi-analytic functions in a domain $\Omega$.
\begin{cor}
	Let $F$ be a bi-analytic function in a domain $\Omega$ with the series expansion as in Theorem \ref{him-P7-thm-1.7}. Also assume all the hypothesis as in Theorem \ref{him-P7-thm-1.7}. Then $M_{r}(F)\leq 1$ for $|z|=r \leq \min \{r_{f}(2), 1/(1+2\lambda _{\mathcal{H}})\}$, where $r_{f}(2)$ is the smallest root in $(0,1)$ of 
	\begin{equation} \label{him-P7-e-1.74}
	(1-r)^2 - k\lambda_{\mathcal{H}}\, r + k\lambda_{\mathcal{H}}\, r^{3}=0,
	\end{equation}
	where $\lambda _{\mathcal{H}}$ is given by \eqref{him-P7-e-1.24}.
\end{cor} 
For $\Omega=\Omega_{\gamma}$, we have $\lambda_{\mathcal{H}}=\lambda_{\mathcal{H}}(\Omega_{\gamma})\leq 1/(1+\gamma)$ (see \cite{Himadri-Vasu-P7}). In view of Theorem \ref{him-P7-thm-1.7}, we obtain the following corollaries.
\begin{cor} \label{him-P7-cor-1.75}
	Let $F$ be a polyanalytic function in $\Omega_{\gamma}$ with the series expansion as in Theorem \ref{him-P7-thm-1.7}. Also assume all the hypothesis as in Theorem \ref{him-P7-thm-1.7}.  Then $M_{r}(F)\leq 1$ for $|z|=r \leq \min \{r_{f}(p,\gamma), (1+\gamma)/(3+\gamma)\}$, where $r_{f}(p,\gamma)$ is the smallest root in $(0,1)$ of 
	\begin{equation} \label{him-P7-e-1.78}
	(1+\gamma)(1-r)^2 -  kr + k r^{p+1}=0.
	\end{equation} 
\end{cor} 
The following result is the limiting case of Corollary \ref{him-P7-cor-1.75}. Consider the domain $\tilde{\Omega}=\{z: \real z<1\}$ which can be obtained as the limiting case of the domain $\Omega_{\gamma}$ by considering $\gamma \rightarrow 1^{-}$.
\begin{cor}\label{him-P7-cor-1.77}
	Let $F$ be a polyanalytic function in $\tilde{\Omega}$ with the series expansion as in Theorem \ref{him-P7-thm-1.7}. Also assume all the hypothesis as in Theorem \ref{him-P7-thm-1.7}.
  Then $M_{r}(F)\leq 1$ for $|z|=r \leq \min \{r_{f}(p,1), 1/2\}$, where $r_{f}(p,1)$ is the smallest root in $(0,1)$ of 
  \begin{equation} \label{him-P7-e-1.76}
  2(1-r)^2 -  kr + k r^{p+1}=0.
  \end{equation}
\end{cor} 
In the next result, we obtain Bohr radius for the polyanalytic function $F(z)= \sum_{l=0}^{p-1} \overline{z}^l \, f_{l}(z) $, where $f_{0}$ is subordinate to a convex biholomorphic function in the unit disk $\mathbb{D}$. 
\begin{thm} \label{him-P8-thm-1.2}
	Let $F$ be a polyanalytic function of order $p$ in $\mathbb{D}$ with $F(z)= \sum_{l=0}^{p-1} \overline{z}^l \, f_{l}(z) $, where each $f_{l}: \mathbb{D} \rightarrow \mathcal{B}(\mathcal{H})$ is analytic function such that $f_{l}(z)= \sum_{n=0}^{\infty} A_{n,l} z^n$ in $\mathbb{D}$ and $A_{n,l} \in \mathcal{B}(\mathcal{H})$ for $n \in \mathbb{N} \cup \{0\}$. Also assume that 
	\begin{enumerate}
		\item $f_{0} \in S(g)$ such that $f_{0}(0)={\bf 0}$ and $f'_{l}(0)=\alpha_{l} f'_{0}(0)$ with $|\alpha_{l}| <k$ for $k \in [0,1]$ and each $l=1,\ldots, p-1$, where $g:\mathbb{D} \rightarrow \mathcal{B}(\mathcal{H})$ is a convex biholomorphic function with $g(z)= \sum_{n=0}^{\infty} g_{n} z^n$ in $\mathbb{D}$ and $g_{n} \in \mathcal{B}(\mathcal{H})$ for $n \in \mathbb{N} \cup \{0\}$.
		\item $\omega_{l} : \mathbb{D} \rightarrow \mathcal{B}(\mathcal{H})$ is analytic function with $\norm{\omega_{l}(z)} \leq k$ in $\mathbb{D}$ for $k \in [0,1]$, where $\omega_{l}(z)= f'_{l}(z) (f'_{0}(z))^{-1}$ in $\mathbb{D}$ such that $\omega_{l}(z)=\sum_{n=0}^{\infty} \omega_{n,l} z^n $ in $\mathbb{D}$, provided $(f'_{0}(z))^{-1}$ exists for all  $z \in \mathbb{D}$.
	\end{enumerate}
	Then $M_{r}(F) \leq 1$ for $|z|=r \leq R_{C}= \min \{r_{C}(p,k,\beta), 1/3\}$, where $r_{C}(p,k,\beta)$ is the smallest root in $(0,1)$ of 
	\begin{equation} \label{him-P8-e-1.14}
	(1-r)^2 -k \beta\, r + k \beta\, r^{p+1}=0,
	\end{equation}
	where $\norm{g'(0)}=\beta$.
\end{thm}
\begin{pf}
From \eqref{him-P7-e-1.70}, it is enough to estimate the upper bound of $M_{r}(f_{0})$. Let $g: \mathbb{D} \rightarrow \mathcal{B}(\mathcal{H})$ be univalent and convex biholomorphic function in $\mathbb{D}$ such that $g(z)=\sum_{n=0}^{\infty} g_{n}z^n$, where $g_{n} \in \mathcal{B}(\mathcal{H})$. Set $\xi = e^{2\pi i/n}$. Since $g$ is convex then by the similar argument used in proving \cite[Theorem X]{rogosinski-1943}, we obtain 
\begin{equation*}
\Psi(z^n)= \frac{f_{0}(\xi z)+f_{0}(\xi^2 z)+ \cdots + f_{0}(\xi^n z)}{n}=A_{n,0}z^n + A_{2n,0}z^{2n}+ \cdots \prec g(z),
\end{equation*}
and hence, $\Psi(z)=A_{n,0}z + A_{2n,0}z^{2}+ \cdots \prec g(z)$ for $z \in \mathbb{D}$. Hence there exits a holomorphic function $\omega : \mathbb{D} \rightarrow \mathbb{D}$ with $\omega(0)=0$ such that $\Psi(z)=g(\omega(z))$, which implies that $\Psi'(0)=\omega'(0)g'(0)$. That is, $A_{n,0}=\omega'(0)g'(0)$, which leads to $\norm{A_{n,0}} \leq \norm{g'(0)}$. Using this inequality and the fact $f_{0}(0)={\bf 0}$, we obtain 
\begin{equation} \label{him-vasu-P6-e-2.51}
M_{r}(f_{0})=\sum_{n=0}^{\infty}\norm{A_{n,0}}r^n \leq \left(\frac{r}{1-r}\right) \, \norm{g'(0)}=\beta \left(\frac{r}{1-r}\right).
\end{equation}
For $\Omega_{\gamma}=\mathbb{D}$ {\it i.e.,} $\gamma=0$ and if $h:\mathbb{D} \rightarrow \mathcal{B}(\mathcal{H})$ is analytic function then we have $\lambda_{\mathcal{H}} \leq 1$ (see \cite{Himadri-Vasu-P7}). Hence by the analogues proof of Theorem \ref{him-P7-thm-1.7}, from \eqref{him-P7-e-1.70}, we obtain 
\begin{equation*}
M_{r}(F) \leq k \beta \left(\frac{r}{1-r}\right) \left(\frac{1-r^p}{1-r}\right)\,\,\,\, \mbox{for} \,\,\,\, r \leq \frac{1}{3}.
\end{equation*}
Hence, $M_{r}(F) \leq 1$ for $r \leq R_{C}= \min \{r_{C}(p,k,\beta), 1/3\}$, where $r_{C}(p,k,\beta)$ is the smallest root in $(0,1)$ of \eqref{him-P8-e-1.14}. This completes the proof.
\end{pf}
\vspace{-4mm}

Let $h:\mathbb{D} \rightarrow \mathcal{B}(\mathcal{H})$ be holomorphic and $g \in S(h)$ with the expansions $h(z)=\sum_{n=0}^{\infty} h_{n}z^n$ and $g(z)=\sum_{n=0}^{\infty} g_{n}z^n$ respectively in $\mathbb{D}$, where $h_{n},g_{n} \in \mathcal{B}(\mathcal{H})$ for $n \in \mathbb{N} \cup \{0\}$. Then, in view of \cite[Lemma 2]{bhowmik-2021}, for $|z|=r\leq 1/3$, it is known that 
\begin{equation} \label{him-P8-e-1.15}
\sum_{n=1}^{\infty} \norm{g_{n}} r^n \leq \sum_{n=1}^{\infty} \norm{h_{n}} r^n.
\end{equation}
In the following result, we obtain Bohr radius for the polyanalytic function $F(z)= \sum_{l=0}^{p-1} \overline{z}^l \, f_{l}(z) $, where $f_{0}$ is subordinate to a starlike biholomorphic function in the unit disk $\mathbb{D}$.
\begin{thm} \label{him-P8-thm-1.3}
	Let $F$ be a polyanalytic function of order $p$ in $\mathbb{D}$ with $F(z)= \sum_{l=0}^{p-1} \overline{z}^l \, f_{l}(z) $, where each $f_{l}: \mathbb{D} \rightarrow \mathcal{B}(\mathcal{H})$ are analytic functions such that $f_{l}(z)= \sum_{n=0}^{\infty} A_{n,l} z^n$ in $\mathbb{D}$ and $A_{n,l} \in \mathcal{B}(\mathcal{H})$ for $n \in \mathbb{N} \cup \{0\}$. Also assume that 
	\begin{enumerate}
		\item $f_{0} \in S(g)$ such that $f_{0}(0)={\bf 0}$ and $f'_{l}(0)=\alpha_{l} f'_{0}(0)$ with $|\alpha_{l}| <k$ for $k \in [0,1]$ and each $l=1,\ldots, p-1$, where $g:\mathbb{D} \rightarrow \mathcal{B}(\mathcal{H})$ is a normalized starlike biholomorphic function with $g(z)= zI+\sum_{n=2}^{\infty} g_{n} z^n$ in $\mathbb{D}$ and $g_{n} \in \mathcal{B}(\mathcal{H})$ for $n \in \mathbb{N} \cup \{0\}$.
		\item $\omega_{l} : \mathbb{D} \rightarrow \mathcal{B}(\mathcal{H})$ is analytic function with $\norm{\omega_{l}(z)} \leq k$ in $\mathbb{D}$ for $k \in [0,1]$, where $\omega_{l}(z)= f'_{l}(z) (f'_{0}(z))^{-1}$ in $\mathbb{D}$ such that $\omega_{l}(z)=\sum_{n=0}^{\infty} \omega_{n,l} z^n $ in $\mathbb{D}$.
	\end{enumerate}
	Then $M_{r}(F) \leq 1$ for $|z|=r \leq R_{S}= \min \{r_{S}(p,k), 1/3\}$, where $r_{s}(p,k)$ is the smallest root in $(0,1)$ of 
	\begin{equation} \label{him-P8-e-1.16}
	(1-r)^3 -k r + k r^{p+1}=0.
	\end{equation}
\end{thm}
\begin{pf}
Let $g:\mathbb{D} \rightarrow \mathcal{B}(\mathcal{H})$ be a normalized starlike biholomorphic function. Then in view of \cite[Theorem 6.2.6]{graham-2003}, $g$ satisfies 
\begin{equation} \label{him-P8-e-1.17}
z\, g'(z)= q(z)g(z)\,\,\,\,\,\,\mbox{for} \,\,\,\, z \in \mathbb{D},
\end{equation}
where $q:\mathbb{D} \rightarrow \mathbb{C}$ is an analytic function with $\real q(z)>0$ in $\mathbb{D}$ and $q(0)=1$. By comparing the coefficients in the power series of both the sides of \eqref{him-P8-e-1.17}, we obtain 
\begin{equation} \label{him-P8-e-1.18}
(n-1)g_{n} = g_{n-1}q_{1} + g_{n-2}q_{2}+\ldots + q_{n-1}\,\,\,\,\,\, \mbox{for} \,\,\,\, n\geq 2.
\end{equation}
By using induction and \eqref{him-P8-e-1.18}, we obtain 
\begin{equation*}
(n-1) \norm{g_{n}} \leq 2(n-1+n-2+\ldots+1)I=n(n-1)I,
\end{equation*}
which turns out that $\norm{g_{n}} \leq n$ for all $n \geq 2$. Since, $f_{0} \in S(g)$, by using \eqref{him-P8-e-1.15}, for $r\leq 1/3$, we obtain
 \begin{equation}
 M_{r}(f_{0})=\sum_{n=1}^{\infty} \norm{A_{n,0}} r^n \leq \sum_{n=1}^{\infty} \norm{g_{n}} r^n \leq r+ \sum_{n=2}^{\infty} nr^n=\frac{r}{(1-r)^2}.
 \end{equation}
 It is known that, for $\Omega_{\gamma}=\mathbb{D}$ {\it i.e.,} $\gamma=0$, if $h:\mathbb{D} \rightarrow \mathcal{B}(\mathcal{H})$ is holomorphic then we have $\lambda_{\mathcal{H}} \leq 1$ (see \cite{Himadri-Vasu-P7}). Then by the analogues proof of Theorem \ref{him-P7-thm-1.7}, from \eqref{him-P7-e-1.70}, we obtain 
 \begin{equation*}
 M_{r}(F) \leq k\,  \frac{r}{(1-r)^2} \left(\frac{1-r^p}{1-r}\right)\,\,\,\, \mbox{for} \,\,\,\, r \leq \frac{1}{3}.
 \end{equation*}
 Hence, $M_{r}(F) \leq 1$ for $r \leq R_{S}= \min \{r_{S}(p,k), 1/3\}$, where $r_{C}(p,k,\beta)$ is the smallest root in $(0,1)$ of \eqref{him-P8-e-1.16}. This completes the proof.
\end{pf}

\noindent\textbf{Acknowledgment:} 
% The authors would like to express their sincerest gratitude to the referee for careful reading of the manuscript and many valuable suggestions, which greatly helped to improve the clarity of the exposition in this manuscript. 
The first author is supported by SERB-CRG and the second author is supported by CSIR (File No: 09/1059(0020)/2018-EMR-I), New Delhi, India.

\end{document}